\documentclass[12pt]{article}
 \usepackage{amsfonts,amssymb,amsmath,mathrsfs}
 \usepackage{color}
 \usepackage[svgnames]{xcolor}

 \newtheorem{theorem}{Theorem}[section]
 \newtheorem{lemma}[theorem]{Lemma}
 \newtheorem{corol}[theorem]{Corollary}
 \newtheorem{prop}[theorem]{Proposition}
 \newtheorem{remark}[theorem]{Remark}
 \newtheorem{example}[theorem]{Example}
 \newtheorem{condition}[theorem]{Condition}

 \def\btheorem{\begin{theorem}\sl{}
 \def\etheorem{\end{theorem}}}
 \def\blemma{\begin{lemma}\sl{}
 \def\elemma{\end{lemma}}}
 
 \def\bproposition{\begin{prop}\sl{}
 \def\eproposition{\end{prop}}}
 \def\bremark{\begin{remark}\sl{}
 \def\eremark{\end{remark}}}

 \def\benumerate{\begin{enumerate}}\def\eenumerate{\end{enumerate}}
 \def\bitemize{\begin{itemize}}\def\eitemize{\end{itemize}}
 \def\itm{\item}

 \def\beqlb{\begin{eqnarray}}
 \def\eeqlb{\end{eqnarray}}
 \def\beqnn{\begin{eqnarray*}}
 \def\eeqnn{\end{eqnarray*}}

 \def\pf{\noindent{\it Proof.~~}}\def\proof{\noindent{\it Proof.~~}}
 \def\qed{\hfill$\Box$\medskip}

 \def\<{\langle}\def\>{\rangle}

 \def\mbb{\mathbb}
 \def\mbf{\mathbf}\def\mrm{\mathrm}

 \def\ar{\!\!&}\def\nnm{\nonumber}\def\ccr{\nnm\\}

 \def\d{\mrm{d}}\def\e{\mrm{e}}

\begin{document}


\bigskip\bigskip

\centerline{\Large\bf A scaling limit theorem for Galton-Watson}

\smallskip

\centerline{\Large\bf processes in varying environments\,\footnote{Supported by the National Key R{\&}D Program of China (No.~2020YFA0712900), the National Natural Science Foundation of China (No.~11531001), the Program for Probability and Statistics: Theory and Application (No.~IRTL1704), the Program for Innovative Research Team in Science and Technology in Fujian Province University (IRTSTFJ) and the Education and Scientific Research Project for Young and Middle-aged Teachers in Fujian Province of China (No.~JAT200072).}}

\bigskip

\centerline{Rongjuan Fang$^{\rm(a)}$, Zenghu Li$^{\rm(b)}$ and Jiawei Liu$^{\rm(c)}$\footnote{Corresponding author.}}

\bigskip

\centerline{(a) College of Mathematics and Informatics,}

\centerline{Fujian Normal University, Fuzhou 350007, P.R.\ China}
\centerline{(b,c) Laboratory of Mathematics and Complex Systems (MOE),}

\centerline{School of Mathematical Sciences,}

\centerline{Beijing Normal University, Beijing 100875, P.R.\ China}

\bigskip

\centerline{E-mails: {\tt fangrj@fjnu.edu.cn}, {\tt lizh@bnu.edu.cn} and {\tt jwliu@mail.bnu.edu.cn}}

\bigskip\bigskip

{\narrower{\narrower

\noindent\textit{Abstract}: {We prove a scaling limit theorem for discrete Galton-Watson processes in varying environments. A simple sufficient condition for the weak convergence in the Skorokhod space is given in terms of probability generating functions. The limit theorem gives rise to the continuous-state branching processes in varying environments studied recently by several authors.}

\medskip

\noindent\textit{Key words and phrases}: Galton-Watson processes; continuous-state; varying environments;  probability generating functions; scaling limits.

\par}\par}

\bigskip\bigskip

\section{Introduction}

\setcounter{equation}{0}

Let $\{g_{n}:n\geq0\}$ be a sequence of probability generating functions. {  Suppose that $\{\xi^{(i)}_{n}:n=0,1,2\cdots  ; i=1,2,\cdots\}$
are indenpendent positive integer-valued  random variables, where $\xi^{(i)}_{n}$ has distribution given by $g_{n}$.}
A {\it Galton-Watson process in varying environments} (GWVE process) $\{Z(n):n\geq0\}$ corresponding to $\{g_{n}\}$ is defined inductively as follows:
\beqnn
Z(n+1)=\sum_{i=1}^{Z(n)}\xi^{(i)}_{n},~~n=0,1,2,\ldots,
\eeqnn
where $\xi^{(i)}_{n}$ means the number of offspring generated by the $i$-th individual in the $n$-th generation.
Therefore $\{Z(n):n\geq0\}$ is a discrete-time non-negative integer-valued inhomogeneous Markov chain with
transition matrix $P(n,i;n+1,j)$ defined by
\beqlb\label{Pk1}
 \sum_{j=0}^{\infty}P(n,i;n+1,j)z^{j}=(g_{n}(z))^{i},~~n,i=0,1,2,\ldots,~|z|\leq 1,
\eeqlb
and  with the multistep transition matrix $P(m,i;n,j)$ determined by
\beqlb\label{Pkmn}
 \sum_{j=0}^{\infty}P(m,i;n,j)z^{j}=(g_{m,n}(z))^{i},~~n\geq m\ge0,~i=0,1,2,\ldots,~|z|\leq 1,
\eeqlb
where $$g_{m,n}(z):=g_{m}(g_{m+1}(\cdots(g_{n-1}(z))\cdots))$$ if $n\geq m+1$ and $g_{m,m}(z):=z$.
If $g_{n}\equiv g$, then $\{Z(n):n\geq0\}$ is the classical {\it Galton-Watson process} (GW process). The research
on GWVE processes dates back to the 1970s. Church (1971) firstly investigated the limiting distribution for the GWVE processes
from a purely analytic perspective. The problem related to the extinction probability was studied in Agresti (1975), Fujimagari (1980)
and others. We refer to Lindvall (1974), Macphee and Schuh (1983) and others for more limit theorems.
The study of GWVE processes is conducive to GW processes in random environments; see, e.g.,
Afanasyev et al. (2005, 2012), Bansaye and B\"{o}inghoff (2013), Kersting and Vatutin (2017) for more details.

The scaling limit theorems of discrete-state processes have attracted widespread attention. The convergence of the rescaled GW processes
was pioneered by Feller (1951). It is well-known that {\it continuous-state branching processes} (CB processes) arise in limit theorems
of rescaled GW-processes; see Aliev and Shchurenkov (1982), Grimvall (1974), Lamperti (1967) and Li (2006). The basic structures of general
CB-processes were studied by Ji\v{r}ina (1958). CB processes without or with immigration can be constructed by stochastic
equations; see, e.g. Dawson and Li (2006, 2012). In regard to GWVE processes, Bansaye and Simatos (2015) unfolded a general sufficient
condition for the weak convergence of the rescaled sequence and characterized the limit process by the cumulant semigroup. They provided
an exhaustive result on the offspring distributions with infinite variances resulting in diffusion processes with jumps. Their condition
is intuitively relying on a triplet of two real functions and a measure but not easy to verify and hence the existence of the limit process
is not obvious. Such a limit process is an inhomogeneous Markov process with branching property, called
a \textit{CB process in varying environments} (CBVE process). Fang and Li (2022+) gave a direct construction of the cumulant
semigroup and the transition semigroup of the CBVE process and established the stochastic integral equation for the process,
which also verified the existence of the limit process given in Bansaye and Simatos (2015). However,  the relationship between
GWVE processes and CBVE processes has not completely been solved yet.  In view of this, we will focus on establishing a more explicit connection.

In this paper, we provide a sufficient condition for the weak convergence of rescaled GWVE processes to the CBVE process in terms of
probability generating functions and illustrate how such conditions can be realized. The remainder of this paper is organized as follows.
In Section~2, we recall some preliminaries about the CBVE process and give the main limit theorem with the sufficient condition. In Section~3,
we find appropriate rescaled GWVE processes satisfying the condition by constructing probability generating functions. In Section~4, we give
a decomposition of the discrete cumulant semigroup and consider the approximation of different parts in the  decomposition. The proofs of the
main theorems are presented in Section~5.

\section{Preliminaries and main results}

\setcounter{equation}{0}

In this section, we first recall some preliminaries about the CBVE-process and then shed light on the scaling procedure to give the main results.

Given a c\`{a}dl\`{a}g function $f$ on $[0,\infty)$ with locally bounded variations and $f(0)= 0$, we write $\Delta f(t)= f(t) - f(t-)$ for the
size of its jump at $t>0$ and $\|f\|(t)$ for the total variation of $f$ on $[0,t]$. Denote $f_c(t)$ be the continuous part of  $f$, i.e., $f_c(t)=f(t)-\sum_{s\leq t} \Delta f(s)$. Let $\tilde{b}_1$ and $\tilde{c}$ be c\`{a}dl\`{a}g functions
on $[0,\infty)$ satisfying $\tilde{b}_1(0)= \tilde{c}(0)= 0$ and having locally bounded variations. Let $\tilde{m}$ be a $\sigma$-finite measure
on $(0,\infty)^2$ satisfying
 \beqlb\label{m.fm}
\tilde{m}_1(t)=\int_0^t\int_0^\infty (1\wedge z^2) \tilde{m}(\d s,\d z)< \infty, \qquad t\ge 0.
 \eeqlb
Here and in the sequel, we understand, for $t\ge r\ge0$,
 \beqnn
\int_r^t =  \int_{(r,t]},\quad \int_r^{t-} =  \int_{(r,t)}
 \quad\mbox{and}\quad
\int_r^\infty  = \int_{(r,\infty)}.
 \eeqnn
{  We introduce the concept of {\it admissible parameters} $(\tilde{b}_1,\tilde{c},\tilde{m})$ defined in Fang and Li (2022+):}

{\rm(1)} $t\mapsto \tilde{c}(t)$ is increasing and continuous;

{\rm(2)} for every $s>0$ we have
 \beqlb\label{Delta0}
\delta(s):= \Delta \tilde{b}_1(s) + \int_0^1 z \tilde{m}(\{s\},\d z)\le 1,
 \eeqlb

{\rm(3)} {  for every $s>0$ we have $ \tilde{m}(\{s\}\times(1,\infty))>0$ when $\Delta \tilde{b}_1(s)=1$.}

{  Consider  the backward integral equation
 \beqlb\label{v.fl}
v(r,t;\lambda)\ar=\ar \lambda - \int_r^t v(s,t;\lambda)\tilde{b}_1(\d s) - \int_r^t v(s,t;\lambda)^2\tilde{c}(\d s) \cr
\ar\ar\qquad - \int_r^t\int_0^\infty K_1(v(s,t;\lambda),z)\tilde{ m}(\d s,\d z),\qquad\lambda>0,
 \eeqlb
where $$K_1(\lambda,z)= \e^{-\lambda z} - 1 + \lambda z \mathbf{1}_{\{z\le 1\}}.$$
Fang and Li (2022+) have proved that there exists a unique bounded positive solution to \eqref{v.fl} with admissible parameters.
Let $\{X(t):t\ge0\}$ be an inhomogeneous Markov process with transition semigroup $\{Q_{r,t}: t\ge r\ge0\}$ defined by
 \beqlb\label{Qrt}
\int_{[0,\infty]}\e^{-\lambda y}Q_{r,t}(x,\d y)= \e^{-xv (r,t;\lambda)}, \qquad \lambda> 0, x\in [0,\infty],
 \eeqlb
with $\e^{-\lambda y}= 0$ for $y=\infty$ by convention and $\{v(r,t): t\ge r\ge0\}$, called the {\it cumulant semigroup},
 defined by the solution to \eqref{v.fl}. The existence of process $\{X_t\}$ is  shown in Fang and Li (2022+).}

In order to construct rescaled discrete processes conveniently, we write \eqref{v.fl} with a common time
scale $\gamma$. Let $s\mapsto\gamma(s)$ be an increasing {  c\`{a}dl\`{a}g} function
on $[0,\infty)$ with $\gamma(0)=0$. Then $\gamma$ induces a measure on $[0,\infty)$.  Suppose that $b_1$ and $c$ are {  bounded}
Borel function on $\mbb{R}_{+}$, and $(1\wedge z^2)m(s,\d z)$ is a {  bounded}
kernel from $\mbb{R}_{+}$ to $\mathscr{B}((0,\infty))$. A time-dependent function $\phi(s,\lambda)$ is given by
\beqnn
\phi(s,\lambda) = b_1(s) \lambda + c(s)\lambda^2 + \int_0^{\infty}K_1(\lambda,z)m(s,\d z).
\eeqnn
In accordance with the admissible condition as in Fang and Li (2022+), $ (b_1,c,m)$ satisfy the following conditions:

{\rm($1'$)} $c$ is non-negative and $c(s)\gamma(\{s\})=0$ for $s>0$;

{\rm($2'$)} for every $s>0$ we have
 \beqlb\label{Delta}
\delta(s):=\Big[ b_1(s)+ \int_0^1 z m(s,\d z)\Big]\gamma(\{s\})\le 1.
 \eeqlb

{\rm($3'$)} {  for every $s>0$ we have $  m(s,(1,\infty))\gamma(\{s\})>0$ when $b_1(s)\gamma(\{s\})=1$.}

Then we claim that \eqref{v.fl} is equivalent to
\beqlb\label{v.fll}
v(r,t;\lambda)\ar=\ar \lambda - \int_r^t \phi(s,v(s,t;\lambda))\gamma(\d s),\qquad\lambda>0,
 \eeqlb
and hence \eqref{v.fll} has a unique bounded positive solution.
\blemma The equation \eqref{v.fl} is equivalent to \eqref{v.fll}.
\elemma

\pf Consider the equation \eqref{v.fl}. For given $T\ge 0$, without loss of generality, we suppose
that $\tilde{m}_1(T)>0$. Then
\beqnn
M_{T}(\d s,\d z):=\tilde{m}_1(T)^{-1}\mbf{1}_{\{s\leq T\}}(1\wedge z^2) \tilde{m}(\d s,\d z)
\eeqnn
is a probability measure on $[0,T]\times(0,\infty)$ by \eqref{m.fm}. Let $\mu_{T}(\d s) = M_{T}(\d s,(0,\infty))$
be the marginal distribution. Hence there exists a transition  probability $K_T(s,\d z)$, which is the regular
conditional distribution, such that
\beqlb\label{decomposition0}
M_{T}(\d s,\d z)=\mu_{T}(\d s)K_T(s,\d z).
\eeqlb
It is obvious that
\beqlb\label{decomposition1}
\mbf{1}_{\{s\leq T\}}\tilde{m}(\d s,\d z)\ar=\ar\tilde{m}_1(T)(1\wedge z^2)^{-1} M_{T}(\d s,\d z) \ccr
 \ar=\ar
\tilde{m}_1(T)\mu_{T}(\d s)(1\wedge z^2)^{-1}K_T(s,\d z).
\eeqlb
Let
$${m}'_T(\d s)=\tilde{m}_1(T)\mu_T(\d s),\ \quad m_T''(s,\d z)=(1\wedge z^2)^{-1}K_T(s,\d z).$$
Then \eqref{decomposition1} implies that
\beqlb \label{decomposition2}
\mbf{1}_{\{s\leq T\}}\tilde{m}(\d s,\d z)={m}'_T(\d s)m_T''(s,\d z).
\eeqlb
By definition, for $0\leq t\leq T$, we have
\beqnn
{m}'_T(0,t]=\tilde{m}_1(T)\mu_T(0,t]=\tilde{m}_1(T)M_T\big((0,t]\times (0,\infty)\big)=\tilde{m}_1(t),
\eeqnn
which implies that ${m}'_T(\d s)$ does not depend on $T$. Hence we can extend the measure ${m}'_T(\d s)$
to the $\sigma$-finite measure $\tilde{m}_1(\d s)$ on $(0,\infty)$ such that ${m}'_T(\d s)=\tilde{m}_1(\d s)\mbf{1}_{\{s\leq T\}}$.
On the other hand,
for $0\leq T_1\leq T_2$, by \eqref{decomposition2},
$$
\mbf{1}_{\{s\leq T_1\}}\tilde{m}(\d s,\d z)=\mbf{1}_{\{s\leq T_1\}}\tilde{m}_1(\d s)m_{T_1}''(s,\d z)
$$
and
$$
\mbf{1}_{\{s\leq T_1\}}\tilde{m}(\d s,\d z)=\mbf{1}_{\{s\leq T_1\}}\mbf{1}_{\{s\leq T_2\}}\tilde{m}(\d s,\d z)
=\mbf{1}_{\{s\leq T_1\}}\tilde{m}_1(\d s)m_{T_2}''(s,\d z).
$$
 Since the decomposition of \eqref{decomposition0} is unique,
we have for $s\leq T_1$, ${m}''_{T_1}(s,\d z)={m}''_{T_2}(s,\d z), \tilde{m}_1$-a.e.. Consequently,
 we can find $m_2(s,\d z)$ such that $\tilde{m}_1$-a.e. $m_2(s,\d z)={m}''_{T}(s,\d z)$ for $s\leq T$
 and $\tilde{m}(\d s,\d z)=\tilde{m}_1(\d s)m_2(s,\d z)$. Moreover $\tilde{m}_1$-a.e. for $s\leq T$,
 \beqlb\label{m2}
 \int_0^\infty (1\wedge z^2) m_2(s, \d z)=\int_0^{\infty}(1\wedge z^2){m}''_{T}(s,\d z)=\int_0^{\infty} K_T(s,dz)=1.
  \eeqlb
Let $$\gamma(\d s)= |\tilde{b}_1|(\d s)+ \tilde{c}(\d s)+ \tilde{m}_1(\d s).$$
Then $\gamma(t)=\int_0^t \gamma(\d s)$  is an increasing function on $[0,\infty)$ with $\gamma(0)=0$.
It is obvious that $\tilde{b}_1,\tilde{c},\tilde{m}_1$ are all absolutely continuous with respect to $\gamma$. Therefore
there exist $b,c,m_1$ such that
\beqnn
 \tilde{b}_1(\d s) = b_1(s)\gamma(\d s);\ \quad \tilde{c}(\d s) = c (s) \gamma(\d s);\ \quad \tilde{m}_1(\d s) = m_1(s) \gamma(\d s).
\eeqnn
Moreover, $\gamma$ a.e. $|b_1(s)|\leq 1,~0\leq c(s),m_1(s)\leq 1$.
Let $m(s,\d z)= m_1(s)m_2(s,\d z)$, which is a $\sigma$-finite kernel. Then $\tilde{m}(\d s,\d z)=m(s,\d z)\gamma(\d s)$ and $\gamma$-a.e.
\beqnn
\int_0^{\infty}(1\wedge z^2) m(s,\d z)=\int_0^{\infty}(1\wedge z^2)m_1(s) m_2(s,\d z)
\leq 1
\eeqnn
by \eqref{m2}. Hence we obtain the equivalence.\qed

Next we consider the weak convergence of rescaled GWVE processes. Let $\{g_{k,n}:n\geq0\}$ be a sequence of probability generating
functions for each $k\ge1$. Consider a sequence of GWVE processes $\{Z_k(n):n\geq0\}_{k\ge1}$ corresponding to $\{g_{k,n}\}$.
Then $\{k^{-1}Z_k(n):n\geq0\}$ is an inhomogeneous Markov chain with state space $E_{k}:=\{0,k^{-1},2k^{-1},\ldots\}$ and the multistep
transition probability $Q_{k}(m,x;n,\mathrm{d}y)$ determined by
\beqlb\label{Qkmn}
 \int_{E_{k}}\e^{-\lambda y}Q_{k}(m,x;n,\mathrm{d}y) = (g_{k,m,n}(\e^{- \frac{\lambda }{k}}))^{kx},\qquad\lambda> 0.
\eeqlb
where $$g_{k,m,n}(z):=g_{k,m}(g_{k,m+1}(\cdots(g_{k,n-1}(z))\cdots))$$ if $n\geq m+1$ and $g_{k,m,m}(z):=z$. Let $\{\beta_k\}$ be a
sequence of increasing constants which tend to infinity as $k\rightarrow \infty$. The time scale is defined
by $\gamma_k(t)=\lfloor \beta_k \gamma(t)\rfloor$ for $k\ge 1$ and $t\ge0$.
{  The rescaled discrete processes $\{X_k(s): s\geq 0\}_{k\ge1}$ are given by
\beqnn
X_k(s) = \frac{1}{k}Z_k(\gamma_k(s)),\qquad k\ge1,~s\geq 0.
\eeqnn}
Denote the discrete cumulant semigroup $  \{v_k(r,t):0\le r\le t\}$ by
\beqlb\label{vk.d}
v_k(r,t;\lambda) = -\log\mbf{E}[\e^{-\lambda X_k(t)}|X_k(r)=1] = -k\log g_{k,\gamma_k(r),\gamma_k(t)}(\e^{- \frac{\lambda }{k}}),~\lambda> 0.
\eeqlb
Before we introduce the sufficient condition for the weak convergence of $\{X_k\}$, some notations are given as follows.
Let $\gamma_k^{-1}(i) = \inf\{s\geq 0:\gamma_k(s)\geq i\}$ for $k \ge 1$ and $i \in \mbb{N}$. For $k\ge 1$ and $  0\leq r< t$, let
\beqlb\label{SJ}
\ar\ar J_k(r,t):=\{s\in(r,t]:\Delta\gamma_k(s)\geq1\},\nnm\\
\ar\ar{J}^{+}_k(r,t):=\{s\in(r,t]: \Delta \gamma_k(s)>1\}.
\eeqlb
For $  s\in J_k(r,t)$,  by the definition of $\gamma_k^{-1}$, we have $\gamma_k^{-1}(\gamma_k(s-)+1)=s$ and let
\beqnn
\ar\ar I_{k,s,1}(\lambda) = k\big[g_{k,\gamma_k(s-)}(1-\lambda/k) - (1-\lambda/k)\big] -
\int_{\gamma_k^{-1}(\gamma_k(s-))}^{{s-}} \phi(u,\lambda) \gamma(\d u),\cr\cr
\ar\ar I_{k,s,2}(\lambda) = k\big[g_{k,\gamma_k(s-) + 1,\gamma_k(s)}(1-\lambda/k)
 - (1-\lambda/k)\big] - \phi(s,\lambda)\gamma(\{s\}),\cr\cr
\ar\ar I_{k,s}(\lambda)=|I_{k,s,1}(\lambda)|+|I_{k,s,2}(\lambda)|.
\eeqnn
{  Consider the following condition in terms of probability generating functions:}
\bitemize
 \itm[{\rm(A)}]{  For $M>0$, there exists  a sequence of  non-negative c\`adl\`ag increasing
 functions $\{F_k(s): s\geq 0\}_{k\geq 1}$  such
 that for $s\geq 0$,~$\lim_{k\rightarrow\infty}F_k(s)=0$,
 and for each $T\geq 0$ and an arbitrary sequence of functions $ \{\lambda_k(s):0\leq s\leq T\}\subset[0, M]$,
 \beqnn
 \sum_{s\in J_k(r,t)}I_{k,s}(\lambda_k(s))\leq F_k(t)-F_k(r),\ \quad 0\leq r\leq t\leq T.
 \eeqnn}
\eitemize
In next section we will construct a sequence of rescaled GWVE processes satisfying the condition (A). Under  condition (A), we can obtain
the boundedness  {  and  convergence} of $v_k(r,t;\lambda)$.
\btheorem  \label{prop.vk-M} Suppose that condition (A) holds. Then for $T\geq 0$ and $0<a<b$,
 there exist  $N_1\geq 1$ and $U\ge l>0$ such that
\beqnn
l \le v_k(r,t;\lambda)\le U,\qquad 0\le r\le t\le T,~a\leq\lambda\leq b,~k\geq N_1.
\eeqnn
\etheorem
\btheorem \label{t.vk.v}
Suppose that condition (A) holds. Then for  $T\geq0$ and $0<a<b$, $v_k(r,t;\lambda)$ converges to $v(r,t;\lambda)$ uniformly for
$0\le r\le t\le T$ and $a\leq \lambda\leq b$,  where $v(r,t;\lambda)$ satisfies \eqref{v.fll}.
\etheorem
\bremark If $\gamma(s)$ is a continuous function, then $J^{+}_k(r,t)\equiv\emptyset$ and $I_{k,s,2}(\lambda)\equiv 0$. {  If we further assume that
 $\phi(s,\lambda)$  is continuous with respect to $(s,\lambda)$, then (A) is implied by
\beqnn
\lim_{k\rightarrow\infty}\sup_{s\in[0,T]}\sup_{0\le\lambda\le M}|\phi_k(s,\lambda)-\phi(s,\lambda)|=0,\ \quad T \geq 0,
\eeqnn
where $$\phi_k(s,\lambda) = \beta_k k\big[g_{k,\gamma_k(s-)}(1-\lambda/k) - (1-\lambda/k)\big].$$}
Then we can obtain the limit theorem following a simple inhomogeneous generalization of Section 3.5 in Li (2011). Comparing it with (A),
 one can see that condition (A) emphasizes the convergence of the integral since $\phi(s,\lambda)$ is no longer continuous in our assumption.
\eremark
Define the metric $d(x,y)=|\e^{-x}-\e^{-y}|$ on the space $[0,\infty]$. Let {  $ D[0,\infty)$ be the space of c\`adl\`ag functions $f:[0,\infty)\mapsto[0,\infty]$}
endowed with the Skorohod topology. Now we give our main limit theorem.

\btheorem\label{t1.4}
Suppose that condition (A)  holds. If $X_k(0)$ converges to $X(0)$ in distribution, then $  \{X_k(s):s \geq 0\}_{k\geq 1}$ is tight
on {  $D[0,\infty)$  and converges} to $  \{X(s):s\geq 0\}$ in distribution on $  D[0,\infty)$ as $k\rightarrow\infty$.
\etheorem

\section{Rescaled GWVE processes}

\setcounter{equation}{0}

In this section we provide the sequence of rescaled GWVE processes satisfying condition (A) by constructing an appropriate time scale and probability
generating functions. By definition,
\beqlb\label{U}
  C_0:=\sup_{s\in[0,\infty)}\Big[ |b_1(s)| + c(s)+\int_0^\infty(1\wedge z^2)m(s,\d z)\Big]<\infty.
\eeqlb
 The next result shows that essentially all CBVE processes defined in Bansaye and Simatos (2015) and Fang and Li (2022+) arise as scaling
limits of GWVE processes in the fashion of Theorem~\ref{t1.4}.

\bproposition\label{prop.gwve}
  For given $\gamma$ and $\phi$ with admissible parameters, there exists a sequence of $(\beta_k,\{g_{k,i}\})$ satisfying condition (A).
\eproposition
\proof   For $k\ge 1$ and $0\leq r< t\leq T$, let
\beqnn
\ar\ar S_k(r,t):=\{i\in\mbb{N}\cap(\gamma_k(r),\gamma_k(t)]:\exists s\in[r,t),\gamma_k(s)=i-1\},\cr\cr
\ar\ar S_k^c(r,t):= \big(\mbb{N}\cap(\gamma_k(r),\gamma_k(t)]\big)\setminus S_k(r,t).
\eeqnn
To construct $g_{k,i}$ for all $i\leq \gamma_k(T)$, we divide $\mbb{N}\cap(0,\gamma_k(T)]$ into $ S_k(0,T)$ and $ S^c_k(0,T)$.
We consider $I_{k,s,1}$ and $I_{k,s,2}$ separately.

\emph{Step~1.}
Let $\beta_k=4C_0(k+1)$ for $k\geq 1$, and
let $c_k= C_0^{-1}k^{1/3}-1$. For sufficiently large $k$,
we have $1\le c_k\le k$. Let
$$
b_k(s)=b_1(s)-\int_1^{c_k} z m(s,\d z).
$$
Then
$$
|b_k(s)|\leq C_0+c_k \int_1^{c_k}  m(s,\d z)\leq(1+c_k)C_0\leq k^{1/3}.
$$
Let
\beqnn
\phi_{0,k}(s,\lambda)=\phi(s,\lambda)-b_k(s)\lambda -\int_{c_k}^{\infty}(\e^{-\lambda z}-1)m(s,\d z)
= c(s)\lambda^2 + \int_0^{c_k}K( \lambda,  z)m(s,\d z),
\eeqnn
where $K(\lambda,z)=\e^{-\lambda z}-1+\lambda z$. { For $s\in J_k(0,T)$, $\gamma_k(s-)+1\in S_k(0,T)$.
For $i\in S_k(0,T)$}, let
 \beqnn
 \hat{g}_{k,i-1}(u) = u + 2k^{-1}\int_{\gamma_k^{-1}(i-1)}^{\gamma_k^{-1}(i)-}
 \phi_{0,k}(s, k(1-u))\gamma(\d s), \qquad u\leq 1,~k\ge 1.
 \eeqnn
To ensure that
\beqnn
\frac{\d^n}{\d z^n}\hat{g}_{k,i}(0)\geq 0,\qquad n\geq 0,
\eeqnn
by elementary calculations, it's sufficient to show that
\beqnn
2\int_{\gamma_k^{-1}(i-1)}^{\gamma_k^{-1}(i)-}
 \Big[2k c(s)+ \int_0^{c_k} z(1-\e^{-kz})m(s,\d z)\Big]\gamma(\d s)
\leq 1.
\eeqnn
Observe that for $  i\in S_k(0,T)$,
 \beqlb\label{gbk}
 \gamma(\gamma_k^{-1}(i)-)\leq \frac{i}{\beta_k} \quad \mbox{and} \quad\gamma(\gamma_k^{-1}(i-1))\geq \frac{i-1}{\beta_k}.
 \eeqlb
 Consequently,
\beqnn
 \ar\ar 2\int_{\gamma_k^{-1}(i-1)}^{\gamma_k^{-1}(i)-}
\Big[2k c(s) + \int_0^{c_k}z(1-\e^{-kz})m(s,\d z)\Big]\gamma(\d s)\cr\cr
\ar\ar\quad
 \leq  2 \int_{\gamma_k^{-1}(i-1)}^{\gamma_k^{-1}(i)-}[2kc(s)+k\int_0^1 z^2 m(s,\d z)+c_k\int_{1}^{\infty}m(s,\d z) ]\gamma(\d s)\cr\cr
 \ar\ar\quad \leq 4kC_0\int_{\gamma_k^{-1}(i-1)}^{\gamma_k^{-1}(i)-}\gamma(\d s)\leq 4kC_0/\beta_k \leq 1,
\eeqnn
which implies that $\hat{g}_{k,i}$ is a probability generating function  for sufficiently large $k$. It is obvious that
 \beqnn
k [\hat{g}_{k,i-1}(1-\lambda/k)-(1-\lambda/k)] =\int_{\gamma_k^{-1}(i-1)}^{\gamma_k^{-1}(i)-}2\phi_{0,k}(s,\lambda)\gamma(\d s),
 \eeqnn
 for $\lambda\in[0,{M}]$  and sufficiently large $k\geq {M}$. For $i\in S_k(0,T)$, let
 \beqnn
  \tilde{g}_{k,i-1}(u)\ar=\ar
 2\Big[\frac{k^{1/3}}{\beta_k}+ \int_{\gamma_k^{-1}(i-1)}^{\gamma_k^{-1}(i)-}b_k(s)\gamma(\d s)\Big]\cr\cr
  \ar\ar\quad
 + \Big[1- 2\int_{\gamma_k^{-1}(i-1)}^{\gamma_k^{-1}(i)-}b_k(s)\gamma(\d s)-4\frac{k^{1/3}}{\beta_k}\Big]u +2\frac{k^{1/3}}{\beta_k} u^2,
\eeqnn
which is also a probability generating function for sufficiently large $k$.
Then for $i\in S_k(0,T)$, define the sequence $\{g_{k,i-1}\}$ by
\beqnn
 g_{k,i-1}:=\frac{1}{2}\hat{g}_{k,i-1}+\frac{1}{2}\tilde{g}_{k,i-1}.
\eeqnn
For an arbitrary sequence of functions $\{\lambda_k(s): 0\leq s\leq T\}\subset[0, {M}]$,
\beqlb\label{Ik1}
\ar\ar\sum_{s\in J_k(r,t)} |I_{k,s,1}(\lambda_k(s))|\cr
\ar\ar= \sum_{s\in J_k(r,t)} \Big|k\Big[{g}_{k,\gamma_k(s-)}\Big(1-\frac{\lambda_k(s)}{k}\Big) - \Big(1-\frac{\lambda_k(s)}{k}\Big)\Big] -
\int_{\gamma_k^{-1}(\gamma_k(s-))}^{{s-}} \phi(u,\lambda_k(s)) \gamma(\d u)\Big|\cr
\ar\ar\le
\sum_{s\in J_k(r,t)}\Big
 |\int_{\gamma_k^{-1}(\gamma_k(s-))}^{s-}\gamma(\d u)\int_{c_k}^{\infty}
 (\e^{-\lambda_{k}(s)z}-1)m(u,\d z)\Big|
\cr
\ar\ar\qquad+\sum_{s\in J_k(r,t)} \Big|k\Big[\tilde{g}_{k,\gamma_k(s-)}\Big(1-\frac{\lambda_k(s)}{k}\Big) - \Big(1-\frac{\lambda_k(s)}{k}\Big)\Big]
-\int_{\gamma_k^{-1}(\gamma_k(s-))}^{{s-}} b_k(s)\lambda_k(s)  \gamma(\d u)\Big|\cr\cr
  \ar\ar\leq \int_{\gamma_k^{-1}(\gamma_k(r))}^{\gamma_k^{-1}(\gamma_k(t))}\gamma(\d u)
  \int_{c_k}^{\infty}m(u,\d z) + \frac{  {M}^2(\gamma_k(t)-\gamma_k(r))}{k^{2/3}\beta_k},
  \eeqlb
where the last inequality follows from
\beqlb\label{sumofS}
\sum_{s\in J_k(r,t)} 1 \le \gamma_k(t)-\gamma_k(r).
\eeqlb

\emph{Step~2.} For $i\in S^c_k(0,T)$, if $i-1\in S^c_k(0,T)$, let $g_{k,i-1}(u)=u$. Otherwise, it implies that $i-2=\gamma_k(s(i)-)$ for some
$s(i)\in {J}^+_k(0,T)$ and let
\beqnn
g_{k,i-1}(u)=\e^{(u-1)[1-\delta_k(\theta,s(i))]} + k^{-1}\int_{k^{-\theta}}^\infty(\e^{-kz(1-u)}-1)m(s(i),\d z)\gamma(\{s(i)\}),
\eeqnn
where $\theta\in(0,1)$ is fixed and
\beqnn
\delta_k(\theta,s):= \Big [b_1(s) + \int_{k^{-\theta}}^1 z m(s,\d z)\Big]\gamma(\{s\})\le \delta(s)\le 1.
\eeqnn
{  Without loss of generality, we assume that $T=N$ for some positive integer $N$ and $(0,T]=\bigcup_{j=1}^{N} (j-1,j].$ For $i\in S_k^c(j-1,j)$ and $i-1\in S_k(j-1,j)$,
we claim that $g_{k,i-1}$ is a probability generating function for some sufficiently large $  k\geq C_1(j)$.} In fact, $g_{k,i-1}(1)=1$ and
\beqnn
g_{k,i-1}(0)\ar=\ar \exp\big\{-[1-\delta_k(\theta,s(i))]\big\} + k^{-1}\int_{k^{-\theta}}^\infty(\e^{-kz}-1)m(s(i),\d z)\gamma(\{s(i)\})\cr
\ar\ge\ar \exp\big\{-[1-\delta_k(\theta,s(i))]\big\} - k^{-1}\int_{k^{-\theta}}^{\infty} m(s(i),\d z)\gamma(\{s(i)\})\cr
\ar\ge \ar \exp\big\{-[1-\delta_k(\theta,s(i))]\big\}- k^{\theta-1}\int_{k^{-\theta}}^1 z m(s(i),\d z)\gamma(\{s(i)\})\cr
\ar\ar\qquad\qquad \qquad\qquad\qquad\qquad\qquad\qquad \qquad -k^{-1}m(s(i),[1,\infty))\gamma\{ s(i) \}
\\
\ar\ge\ar \exp\big\{-[1-\delta_{b_1}(s(i))]\big\}- k^{\theta-1}[1-\delta_{b_1}(s(i))]{  -k^{-1}C_0\gamma(j)},
\eeqnn
where $\delta_{b_1}(s):=b_1(s)\gamma(\{s\})\le \delta_k(\theta,s)$. Hence $g_{k,i-1}(0)\ge 0$ for some sufficiently large $  k\geq C_1(j)$ by {  the
local  boundedness of $\gamma$}. {  For $k< C_1(j)$, let $g_{k,i-1}(z)=z$ for $i\in S_k^c(j-1,j)$ and $i-1\in S_k(j-1,j)$. Note that for} $s\in J^+_k(0,T)$, $\gamma_k(s-)+2\in S^c_k(0,T)$ and $\gamma_k(s-)+1\in S_k(0,T)$. {  Hence,
for} $0\leq r\leq t\leq T$ and $\{\lambda_{k}(s),0\leq s\leq T\}\subset [0, {M}]$, we have
\beqlb\label{I_k.B}
 \ar\ar\sum_{s\in J_k(r,t)} |I_{k,s,2}(\lambda_k(s))|\cr\cr
 \ar\ar\qquad=\sum_{s\in {J}^+_k(r,t)}\big|k[g_{k,\gamma_k(s-)+1}(1-\lambda_{k}(s)/k)-(1-\lambda_{k}(s)/k)]
 - \phi(s,\lambda_{k}(s))\gamma(\{s\})\big|\cr\cr
 \ar\ar\qquad\qquad +\sum_{s\in(r,t]:\Delta \gamma_k(s)=1}|\phi(s,\lambda_{k}(s))\gamma(\{s\})|\cr\cr
 \ar\ar\qquad\le\sum_{s\in {J}^+_k(r,t)}\Big|k K(\lambda_{k}(s)/k,1-\delta_k(\theta,s))
 -\int_0^{k^{-\theta}}K(\lambda_{k}(s),z)m(s,\d z)\gamma(\{s\})\Big|\cr\cr
\ar\ar\qquad\qquad +C_0(1+M)^2\sum_{s\in(r,t]:\Delta \gamma_k(s)=1}\gamma(\{s\})\cr\cr
 \ar\ar\qquad\le\frac{M^2}{2} \sum_{s\in {J}^+_k(r,t)}\frac{( 1 - \delta_k(\theta,s))^2}{k}
 +\frac{M^2}{2}\sum_{s\in {J}^+_k(r,t)}\int_0^{k^{-\theta}}z^2m(s,\d z)\gamma(\{s\})\cr\cr
 \ar\ar\qquad\qquad +C_0(1+M)^2\sum_{s\in(r,t]}\Delta \gamma(s)\mbf{1}_{\{\Delta\gamma(s)\leq 2/\beta_k\}}\cr\cr
 \ar\ar\qquad\le\frac{M^2}{2} \sum_{s\in  {J}^+_k(r,t)}\frac{( 1 - \delta_{b_1}(s))^2}{k}
 +\frac{M^2}{2}\int_r^t \gamma(\d s) \int_0^{k^{-\theta}}z^2m(s,\d z)\cr\cr
 \ar\ar\qquad\qquad +C_0(1+M)^2\sum_{s\in(r,t]}\Delta \gamma(s)\mbf{1}_{\{\Delta\gamma(s)\leq 2/\beta_k\}}.
 \eeqlb
Note that for $s\in {J}^+_k(r,t)$,
\beqnn
\Delta\gamma(s)\ge \frac{\Delta\gamma_k(s)-1}{\beta_k}\geq\frac{1}{\beta_k}=\frac{1}{4C_0(k+1)}.
\eeqnn
Thus $$\mbf{1}_{\{s\in {J}^+_k(r,t)\}}k^{-1}\le
\mbf{1}_{\{s\in (r,t]:\Delta\gamma(s)>0\}}8C_0\Delta\gamma(s).$$
By the dominated convergence theorem, the last term in \eqref{I_k.B} tends to $0$ as $k\rightarrow\infty$ for $r=0$ and $t=T$.
Then \eqref{Ik1} and \eqref{I_k.B} imply that the sequence $\{g_{k,i-1}\}$ satisfies (A).

{  Finally, $(\beta_k, \{g_{k,i}\})$ constructed above can be compatible with different $T$. Let
\beqnn
F_k(t)\ar=\ar \int_0^{\gamma_k^{-1}(\gamma_k(t))}\gamma(\d u)
  \int_{c_k}^{\infty}m(u,\d z) + \frac
  { M^2 \gamma_k(t)}{k^{2/3}\beta_k}+
  \frac{M^2}{2} \sum_{s\in  J^+_k(0,t)}\frac{( 1 - \delta_{b_1}(s))^2}{k}\cr\cr
  \ar\ar\quad
 +\frac{M^2}{2}\int_0^t \gamma(\d s) \int_0^{k^{-\theta}}z^2m(s,\d z)
  +C_0(1+M)^2\sum_{s\in(0,t]}\Delta \gamma(s)\mbf{1}_{\{\Delta\gamma(s)\leq 2/\beta_k\}}.
\eeqnn
Then we complete the proof.}
\qed

\section{Decomposition and approximation of the discrete cumulant semigroup }

\setcounter{equation}{0}

In this section, we first give a decomposition of $v_k(r,t;\lambda)$ corresponding to the continuous structure and jump structure
of \eqref{v.fll}. Then we consider the convergence of different parts of the  decomposition.

Some key functions related with probability generating functions are introduced as follows. {  For $ s> 0$} and $\lambda>0$, let
\beqnn
 \phi^{(1)}_k(s,\lambda)=\lambda+k\log g_{k,\gamma_k(s-)}(\e^{-\lambda/k}),
\eeqnn
and
\beqnn
\phi^{(2)}_k(s,\lambda)=\lambda+k\log g_{k,\gamma_k(s-)+1,\gamma_k(s)}(\e^{-\lambda/k}).
\eeqnn
{  Recall the definition of $J_k(r,t)$ and ${J}^+_k(r,t)$.
Now, we give the decomposition of $v_k(r,t;\lambda)$.}
\blemma\label{le.v-de} For $k\ge 1$, $\lambda>0$ {  and $0\leq r\leq t$},
\beqlb\label{eq.v-de}
  v_k(r,t;\lambda)  =   \lambda - \sum_{s\in J_k(r,t)} L_{k}(s,t;\lambda) - \sum_{s\in {J}^+_k(r,t)}
 \phi^{(2)}_k\big(s,v_k(s,t;\lambda)\big),
\eeqlb
where
\beqnn
 L_{k}(s,t;\lambda):=v_k(s,t;\lambda) - v_k(s-,t;\lambda) - \phi^{(2)}_k \big(s,v_k(s,t;\lambda)\big).
\eeqnn
\elemma
\pf By the definition \eqref{vk.d}, $s\mapsto v_k(s,t;\lambda)$ is a c\`{a}dl\`{a}g piecewise  constant function, then
 \beqnn
v_k(r,t;\lambda)
\ar =\ar
\lambda-\sum_{s\in J_k(r,t) }
[v_k(s,t;\lambda)-v_k (s-,t;\lambda )]\cr
\ar=\ar
\lambda-\sum_{s\in J_k(r,t)}
 \phi^{(2)}_k\big(s,v_k(s,t;\lambda)\big)-\sum_{s\in J_k(r,t) } L_{k}(s,t;\lambda)\cr\cr
 \ar=\ar
\lambda-\sum_{s\in J^+_k(r,t)}
 \phi^{(2)}_k\big(s,v_k(s,t;\lambda)\big) - \sum_{s\in {J}_k(r,t)} L_{k}(s,t;\lambda),
\eeqnn
where the last equality follows from the fact that $\phi^{(2)}_k (s, \lambda)$ is zero if $\Delta \gamma_k(s)=1$.\qed

{  For $0<\lambda\le k$, define}
\beqlb\label{Phi1}
\Phi^{(1)}_k(s,\lambda)=k [g_{k,\gamma_k(s-)}(1-\lambda/k)-(1-\lambda/k)],
\eeqlb
and
\beqlb\label{Phi2}
\Phi^{(2)}_k(s,\lambda)=k [g_{k,\gamma_k(s-)+1,\gamma_k(s)}(1-\lambda/k)-(1-\lambda/k)].
\eeqlb
By a simple calculation, there is a connection between $\Phi^{(i)}_k$ and $\phi^{(i)}_k$,
\beqlb\label{Phi-phi}
\phi^{(i)}_k(s,\lambda)= k\log\Big[1+\frac{\e^{\lambda/k}}{k}\Phi^{(i)}_k\big(s, k(1-\e^{-\lambda/k})\big)\Big],\ \quad i=1,2.
\eeqlb
Let ${M}_{\phi}(\lambda):=  (1+\lambda)^2 C_0$ and
\beqnn
{M}'_{\phi}(\lambda_1,\lambda_2):= {M}_{\phi}(\lambda_2+ \e^{-1}\lambda_1^{-1})~\mbox{for}~ 0<\lambda_1\le\lambda_2.
\eeqnn
 Then
\beqlb\label{phiv-M}
|\phi(s,\lambda)|\leq  {M}_{\phi}(\lambda),
\eeqlb
and for $0<\lambda_1\le\lambda_2$,
\beqlb\label{dphi-M}
|\phi(s,\lambda_1)-\phi(s,\lambda_2)|\leq {M}'_{\phi}(\lambda_1,\lambda_2)|\lambda_1-\lambda_2|.
\eeqlb
Under the condition (A) given in Section 2, some properties of $\Phi^{(1)}_k$ and $\Phi^{(2)}_k$ are given as follows.
\bproposition\label{prop.Phi-k-M}
Suppose that condition (A) holds. {  Then for ${M}>0$,
there exists a sequence of non-negative c\`adl\`ag increasing functions $\{F_{1,k}(t): t\geq 0\}_{k\geq 1}$
such that for each $T\geq 0$,~$F_{1,k}(t)$ converges to  $2{M}_{\phi}({M})\gamma(t)$ uniformly for $0\leq t\leq T$ as $k\rightarrow\infty$,
and for  an arbitrary sequence of functions $\{\lambda_{k}(s):0\leq s\leq T\}\subset[0, {M}]$,
\beqnn
\sum_{s\in J_k(r,t)\atop i=1,2}
 \big|\Phi_k^{(i)}\big(s,\lambda_k(s)\big)\big|\leq F_{1,k}(t)-F_{1,k}(r),\ \quad 0\leq r\leq t\leq T,~k\geq 1.
 \eeqnn}
\eproposition
\pf Note that for $s \in J_k(r,t)$, $\gamma_k(s-)+1\in S_k(r,t)$. For $k\ge 1$, by \eqref{gbk},
we have
\beqnn
\big|\Phi^{(1)}_k(s,\lambda_k(s))\big|
\ar\le\ar
 |I_{k,s,1}(\lambda_k(s)| + {M}_{\phi}( {M})/\beta_k.
\eeqnn
On the other hand,
\beqnn
\big|\Phi^{(2)}_k(s, {\lambda}_{k}(s))\big|
\ar\le\ar |I_{k,s,2}( {\lambda}_k(s))| +  M_{\phi}({M}) \gamma(\{s\}).
\eeqnn
Then it follows from condition (A) and \eqref{sumofS} that
 \beqnn
 \sum_{s\in J_{k}(r,t)\atop i=1,2}
 \big|\Phi^{(i)}_k(s,\lambda_{k}(s))\big|
 \ar\leq\ar
 F_k(t)-F_k(r)+\frac{ {M}_{\phi}({M})}{\beta_k}(\gamma_k(t)-\gamma_k(r))+M_{\phi}({M})\big(\gamma(t)-\gamma(r)\big).
 \eeqnn
Condition (A) implies that $F_k(t)$ tends to $0$ uniformly for $0\leq t\leq T$ as $k\rightarrow\infty$. It's easy to
check that $\beta_k^{-1}\gamma_k(t)$ converges to $\gamma(t)$ uniformly for $t\ge 0$ as $k\rightarrow \infty$, which leads to the desired result with
$$
F_{1,k}(t):= F_k(t)+\frac{ {M}_{\phi}({M})}{\beta_k}\gamma_k(t)+M_{\phi}({M})\gamma(t).
$$
\qed

In the rest of this section, we will assume  the following condition:
\bitemize
\itm[{\rm(B)}] For $T\geq 0$ and $0<a<b$,
 there exist $N_2\geq l^{-1}U^2$ and $U\ge l>0$ such that
\beqnn
l \le v_k(r,t;\lambda)\leq U,\qquad  0\le r\le t\le T,~a\leq\lambda\leq b,~k\geq N_2.
\eeqnn
\eitemize
Suppose that (B) holds. Then for $  k\geq N_2$ and $s\in [r,t]$,
\beqlb\label{vh-M}
l/2\le h_k(v_k(s,t;\lambda))\le v_k(s,t,\lambda) \leq U,
\eeqlb
where
$$z-z^2/2k\le h_k(z):=k(1-\e^{-z/k})\le z.$$
{  Without loss of generality, under (B), we can further assume that
\beqlb\label{bd}
h_k(v_k(s,t;\lambda))\le v_k(s,t;\lambda)\leq U,\ \quad 0\le r\le t\le T,~a\leq\lambda\leq b,~k\geq 1,
\eeqlb
since $v_k(s,t;\lambda)\le v_k(s,t;b)$ by the monotonicity of $v_k$ with respect to $\lambda$ and there are finite terms of $v_k(s,t;b)$ for $k\leq N_2,~0\le s\le t\le T$.
}
Under the conditions (A) and (B), we will prove that when $k$ tends to infinity, the second term of the right-hand
side of \eqref{eq.v-de} will be close to
$$
\sum_{s\in J_k(r,t)} \int_{\gamma_k^{-1}(\gamma_k(s-))}^{s-} \phi\big(s,v_k(s,t;\lambda)\big) \gamma(\d s),
$$
and the last term of that will be close to
$$
\sum_{s\in {J}^+_k(r,t)}\phi\big(s,v_k(s,t;\lambda)\big)\gamma(\{s\}).
$$
For convenience, we give some notations that will be used afterwards. For $0\leq t\leq T,~\lambda> 0$ and $s\in J_k(0,t)$, let
\beqlb\label{vwu}
\ar\ar \tilde{\lambda}_{k}(s,t;\lambda)=v_k(s,t;\lambda);\nnm\\
\ar\ar w_{k}(s,t;\lambda)=-k\log g_{k,\gamma_k(s-)+1,\gamma_k(s)} (\e^{-\tilde{\lambda}_{k}(s,t;\lambda)/k});\nnm\\
\ar\ar  u_{k}(s,t;\lambda)=\e^{-w_{k}(s,t;\lambda)/k}.
\eeqlb
Denote $\tilde{\lambda}_{k}(s,t;\lambda),w_{k}(s,t;\lambda)$ and $u_{k}(s,t;\lambda)$ by $\tilde{\lambda}_{k}(s),w_{k}(s) $ and $u_{k}(s) $ without confusion.
For $0\leq t\leq T$, $\lambda> 0$ and $s\in J_k(0,t)$, we have
\beqnn
\tilde{\lambda}_{k}(s-)\ar=\ar-k\log g_{k,\gamma_k(s-)}(\e^{-w_{k}(s)/k})=
-k\log g_{k,\gamma_k(s-)}
\big(1-k^{-1}h_k(w_{k}(s))\big),
\eeqnn
which implies that for $s\in J_k(0,t)$,
\beqlb\label{w-phi}
h_k(w_{k}(s))-h_k(\tilde{\lambda}_{k}(s-)) =\Phi_k^{(1)}\big(s,h_k(w_{k}(s))\big).
\eeqlb
Next we give the convergence of the jump structure.
\blemma\label{l.cj}
Suppose that conditions (A) and (B) hold. Then
there exists  a sequence of  non-negative c\`adl\`ag increasing functions $\{F_{2,k}(s): s\geq 0\}_{k\geq 1}$
such that for $s\geq 0$,~$\lim_{k\rightarrow\infty}F_{2,k}(s)=0$,
 and for $0\le r\le t\le T,~\lambda\in[a,b]$,
\beqlb\label{eql4.3}
 \sum_{s\in {J}^+_k(r,t)}
 |\phi^{(2)}_k\big(s,v_k(s,t;\lambda)\big)-\phi(s,v_k(s,t;\lambda))\gamma(\{s\})|\leq F_{2,k}(t)-F_{2,k}(r),\ \quad k\geq 1.
 \eeqlb
\elemma
\pf  {  Let $\lambda\in [a,b]$. At the beginning, we assume that $N-1 <r\le t\leq N$
for some integer $N\geq 1$.}
Set $\tilde{\lambda}_k(s)=v_k(s,t;\lambda)$ and $\hat{\lambda}_k(s)=h_k(v_k(s,t;\lambda))$, which satisfy \eqref{vh-M} by (B).
By Proposition \ref{prop.Phi-k-M}, there exists $   {M}_N$  such that,
$$
\sum_{s\in {J}^+_k(r,t)}\big|\Phi^{(2)}_k(s,\hat{\lambda}_{k}(s))\big|\le {M}_N.
$$
{  By condition (B),~\eqref{vh-M} and \eqref{bd}, we can assume that there  exist $U_N,~l_N$, and $C_2(N)$ such that
\beqnn
\ar\ar l_N/2\leq h_k\big(v_k(r,t;\lambda)\big),\ \quad k\geq C_2(N);\cr\ar\ar   v_k(r,t;\lambda)\leq U_N,\ \quad k\geq 1.
\eeqnn
}
Let $H(x)=\log(1+x)-x$. Then by the mean value theorem, for a given $\varepsilon\in(0,1)$ and $|x|\le \varepsilon$,
\beqlb\label{H-2}
|H(x)|\le \frac{x^2}{1-\varepsilon}.
\eeqlb
Given $\varepsilon\in(0,1)$, for some sufficiently large $  k\geq C_2(N)$ and $s\in J^+_k(r,t)$,
$$
\frac{\e^{\tilde{\lambda}_{k}(s)/k}}{k} \Phi^{(2)}_k(s,\hat{\lambda}_{k}(s))\le \varepsilon,
$$
hence by \eqref{Phi-phi} and Proposition \ref{prop.Phi-k-M},
\beqlb\label{1230e1}
\ar\ar\sum_{s\in {J}^+_k(r,t)}\big|\phi^{(2)}_k(s,\tilde{\lambda}_k(s))
- \e^{\tilde{\lambda}_{k}(s)/k}\Phi^{(2)}_k (s, \hat{\lambda}_{k}(s)) )\big|
=   \sum_{s\in {J}^+_k(r,t)} k\Big|H\Big(\frac{\e^{\tilde{\lambda}_{k}(s)/k}}{k} \Phi^{(2)}_k(s,\hat{\lambda}_{k}(s))\Big)\Big|\cr
\ar\ar\qquad\le
\sum_{s\in {J}^+_k(r,t)}\frac{[\e^{\tilde{\lambda}_k(s)/k}\Phi^{(2)}_k (s,\hat{\lambda}_{k}(s) ]^2}{(1-\varepsilon)k}
\le  \frac{\e^{2U_N}M_N}{(1-\varepsilon)k}(F_{1,k}(t)-F_{1,k}(r)).
\eeqlb
Note that for $  k\geq 1$,
\beqlb\label{(1-e)phi}
\sum_{s\in {J}^+_k(r,t)}(\e^{\tilde{\lambda}_{k}(s)/k}-1) \big| \Phi^{(2)}_k (s, \hat{\lambda}_{k}(s))\big|
  \le (\e^{U_N/k}-1)(F_{1,k}(t)-F_{1,k}(r)),
\eeqlb
and by (A),
\beqlb\label{1230e2}
\sum_{s\in {J}^+_k(r,t)}|\Phi^{(2)}_k (s,\hat{\lambda}_{k}(s) ) - \phi (s,\hat{\lambda}_{k}(s) )\gamma(\{s\})|
\leq F_k(t)-F_k(r).
\eeqlb
Moreover, \eqref{dphi-M} implies that {  for $ k\ge C_2(N)$,
\beqlb\label{eq619}
  \ar\ar\sum_{s\in {J}^+_k(r,t)}|\phi(s,\hat{\lambda}_{k}(s)) - \phi(s,\tilde{\lambda}_{k}(s))|\gamma(\{s\})\cr
 \ar\ar\qquad\leq
{M}'_\phi(l_N/2,U_N) \sum_{s\in {J}^+_k(r,t)}|h_k(\tilde{\lambda}_{k}(s))- \tilde{\lambda}_{k}(s)|\gamma(\{s\})\cr
\ar\ar\qquad\leq  {M}'_\phi(l_N/2,U_N) \sum_{s\in {J}^+_k(r,t)}\frac{\tilde{\lambda}_{k}(s)^2}{2k}\gamma(\{s\})
\le \frac{{M}'_\phi(l_N/2,U_N) U_N^2}{2k}\big(\gamma(t)-\gamma(r)\big).
\eeqlb
 Let $F_{2,k}(0)=0$ and  assume that $F_{2,k}(s)$  are well defined on $0\leq N-2<s\leq N-1$  and satisfy \eqref{eql4.3} for $ N-2<r\leq t\leq N-1$.
Now, consider $s\in(N-1,N]$. For $k<C_2(N)$, set
 \beqnn
 F_{2,k}(s)=F_{2,k}(N-1)+\sum_{u\in {J}^+_k(N,s)}
 \sup_{\lambda\leq U_N}|\phi^{(2)}_k(u,\lambda)-\phi(u,\lambda)\gamma(\{u\})|
 \eeqnn
 and for $k\geq C_2(N)$, set
 \beqnn
  F_{2,k}(s)\ar=\ar F_{2,k}(N-1)+\frac{\e^{2U_{N}}M_{N}}{(1-\varepsilon)k}(F_{1,k}(s)-F_{1,k}(N-1))
+(\e^{U_N/k}-1)(F_{1,k}(s)-F_{1,k}(N-1))\cr\cr\ar\ar\quad+F_k(s)-F_k(N-1)
+\frac{{M}'_\phi(l_N/2,U_N) U_N^2}{2k}\big(\gamma(s)-\gamma(N-1)\big).
 \eeqnn
As a result of  \eqref{eq619},~\eqref{1230e1},~\eqref{(1-e)phi} and \eqref{1230e2}, by induction,
we can construct $F_{2,k}(s)$ on $[0,\infty)$ to  obtain the desired result.}
\qed

To obtain the strictly positive lower bound of $|k(1-u_{k}(s,t;\lambda))|$, we can choose a sufficiently large
constant $\eta_T>1$ so that for  $s\in [0,T]$, $m(s, (1,\eta_T])\gamma(\{s\})> 0$ when $b_1(s)\gamma(\{s\})= 1$ by
the admissible condition. Similar to the proof of Proposition 2.3 in Fang and Li (2022+), we introduce a c\`{a}dl\`{a}g
function $ \alpha(r)= \alpha(r,T,U)$ with $\Delta \alpha(s)>-1$:
\beqlb\label{eq2.20}
\alpha(r)\ar=\ar - \frac{1}{2}U \int_0^r\gamma(\d s)\int_0^{\epsilon } z^2 m(s,\d z) + H \int_0^r\gamma(\d s)\int_1^{\eta_T} z m(s,\d z)\cr
 \ar\ar
-\int_0^r b_1(s)\gamma(\d s) - U \int_0^r c(s)\gamma(\d s)
- (1-F)\int_0^r\gamma(\d s)\int_{\epsilon }^1 z  m(s,\d z),
 \eeqlb
and hence
\beqnn
\alpha(\d s)\ar=\ar \Lambda(s)\gamma(\d s),
\eeqnn
where
\beqlb\label{FHe}
F = U ^{-1}(1-\e^{-U }),~
H = (\eta_T U)^{-1}(1-\e^{-\eta_T U }),~
\epsilon = 1\land (U ^{-1}F ),
\eeqlb
and
\beqlb\label{Lambda}
\Lambda(s)=- \frac{1}{2}U\int_0^{\epsilon } z^2 m(s,\d z)+ H\int_1^{\eta_T} z m(s,\d z)-b_1(s) - U  c(s)- (1-F)\int_{\epsilon }^1 z  m(s,\d z).\nnm\\
\eeqlb
\blemma \label{l.u}
Suppose that (A) and (B) hold. Then for $T\geq 0$ and $0<a\leq b $, there exist $M_u\ge l_u>0$
and $N_3\geq N_2\vee l_u^{-1}M_u^2$ such that for  $\lambda\in[a,b]$ and $r\le t\in[0,T]$,
\beqnn
\ar\ar\sup_{k\ge N_3,s\in J_k(r,t)} |k(1-u_{k}(s))|\le M_u,\\
\ar\ar\inf_{k\ge N_3,s\in J_k(r,t)}|k(1-u_{k}(s))|\ge l_u,
\eeqnn
where
\beqlb\label{Mulul'}
 \ar\ar  l_u :=2^{-1}l\prod_{s\in(0,T]} [1+(0\land\Delta\alpha(s))]>0,\cr\cr
 \ar\ar M_u := U + M_{\phi}(U)\gamma(T) + l_u/2,
\eeqlb
and $\alpha$ is defined in \eqref{eq2.20}.
\elemma
\pf Recall that $\tilde{\lambda}_{k}(s)=v_k(s,t;\lambda)$ and $u_{k}(s)= g_{k,\gamma_k(s-)+1,\gamma_k(s)} (\e^{-\tilde{\lambda}_{k}(s)/k})$
for $s\in J_k(r,t)$. If $\Delta\gamma_k(s)=1$, then $u_k(s)=\e^{-v_k(s,t;\lambda)/k}$ and the result  is trivial by \eqref{vh-M}.
Otherwise, it implies that $s\in {J}^+_k(r,t)$ and
\beqnn
\ar\ar u_{k}(s) = \exp\bigl\{-k^{-1}\big[\tilde{\lambda}_{k}(s) - \phi^{(2)}_k(s,\tilde{\lambda}_{k}(s))\big]\bigr\}.
\eeqnn
Note that
\beqnn
 \tilde{\lambda}_{k}(s) - \phi^{(2)}_k(s,\tilde{\lambda}_{k}(s))
 \ar =\ar
 \tilde{\lambda}_{k}(s) -\phi(s,\tilde{\lambda}_{k}(s)) \Delta\gamma(s)  + \phi(s,\tilde{\lambda}_{k}(s))\Delta\gamma(s)
- \phi^{(2)}_k(s,\tilde{\lambda}_{k}(s)).
\eeqnn
It is easy to verify that
\beqnn
 2 l_u\ar \le \ar \tilde{\lambda}_{k}(s) (1+\Delta\alpha(s))
 \le
 \tilde{\lambda}_{k}(s) - \phi(s,\tilde{\lambda}_{k}(s))\Delta\gamma(s)  \le U + M_\phi(U)\gamma(T).
 \eeqnn  By Lemma \ref{l.cj}, we {  have
\beqnn
|\phi(s,\tilde{\lambda}_{k}(s))\Delta\gamma(s)
- \phi^{(2)}_k(s,\tilde{\lambda}_{k}(s))|\le  {F}_{2,k}(T)\le {  l_u}/2,\ \quad k\geq C_3,
\eeqnn
for some $C_3$.}
Then by choosing $  N_3\ge C_3\vee l_u^{-1}M_u^2$, it is elementary to obtain the result.\qed

Now we give the convergence of the continuous structure. Denote $L_{k}(s,t;\lambda)$ by $L_{k}(s)$ without confusion. Recall that
$\tilde{\lambda}_{k}(s)=v_k(s,t;\lambda)$ and $k(1-u_{k}(s))=h_k(w_{k}(s))$.
For $s\in J_k(r,t)$,
\beqlb\label{1230e4}
L_{k}(s) \ar = \ar \tilde{\lambda}_{k}(s) -  \phi^{(2)}_k \big(s,\tilde{\lambda}_{k}(s)\big) +
k\log g_{k,\gamma_k(s-),\gamma_k(t)}(\e^{-\lambda/k})\nnm\\
 \ar = \ar k\log g_{k,\gamma_k(s-)} \circ g_{k,\gamma_k(s-)+1,\gamma_k(s)}(\e^{-\tilde{\lambda}_{k}(s)/k}) - k\log g_{k,\gamma_k(s-)+ 1,\gamma_k(s)}(\e^{-\tilde{\lambda}_{k}(s)/k})\nnm\\
\ar=\ar
k\log \frac{g_{k,\gamma_k(s-)}(u_{k}(s))}{u_{k}(s)} \cr
\ar=\ar
k\log\Big(\frac{\Phi_k^{(1)}\big(s,k(1-u_{k}(s))\big)} {ku_{k}(s)}+1\Big)=k\log\Big(\frac{\Phi_k^{(1)}\big(s,h_k(w_{k}(s))\big)}{ku_{k}(s)}+1\Big).\nnm\\
\eeqlb

\blemma\label{l.cc}
Suppose that conditions (A) and (B) hold. {   Then there exists  a sequence of non-negative c\`adl\`ag increasing
functions $\{{F}_{3,k}(s): s\ge0\}_{k\geq 1}$  such that for $s\geq 0$,~$\lim_{k\rightarrow\infty}F_{3,k}(s)=0$,
 and for $0\le r\le t\le T$ and $\lambda\in[a,b]$,
\beqlb\label{eq1.cc}
 \sum_{s\in J_k(r,t)}
 \Big|L_{k}(s,t;\lambda) - \int_{\gamma_k^{-1}(\gamma_k(s-))}
 ^{s-}\phi(s,v_k(s,t;\lambda))\gamma(\d s)\Big|
 \leq F_{3,k}(t)-F_{3,k}(r),\ \quad k\geq 1.
 \eeqlb}
\elemma
\pf
Recall that $\tilde{\lambda}_{k}(s),~w_{k}(s)$ and  $u_{k}(s)$ are defined in \eqref{vwu}.
{  It's elementary to show that
\beqnn
L_k(s,t;\lambda)=\phi_k^{(1)}\big(s,\tilde{\lambda}_k(s)-\phi_k^{(2)}(s,\tilde{\lambda}_k(s))\big).
\eeqnn
Therefore by repeating the inductive procedures in the proof of Lemma \ref{l.cj}, it's
sufficient to consider the case $N-1<r\leq t\leq N$ and obtain \eqref{eq1.cc} for sufficiently large $k$.
By} Lemma \ref{l.u} and Proposition \ref{prop.Phi-k-M}, there exists $  C_4=C_4(N)$ such
that for sufficently large $k$,~$u_{k}(s)\ge 1/2$,
\beqnn
  \Big|\frac{\Phi_k^{(1)}\big(s, h_k(w_{k}(s))\big)}{ku_{k}(s)}\Big|\le \frac{2C_4}{k}\le \frac{1}{2},
\eeqnn
and hence it follows from \eqref{H-2} and \eqref{1230e4} that
\beqlb\label{1230e5}
\ar\ar\sum_{s\in J_k(r,t)}\Big|L_{k}(s)-\frac{\Phi_k^{(1)}\big(s, h_k(w_{k}(s))\big)}{u_{k}(s)}\Big|
 = \sum_{s\in J_k(r,t)} k \Big|H\Big(\frac{\Phi_k^{(1)}\big(s, h_k(w_{k}(s))\big)}{ku_{k}(s)}\Big) \Big| \cr
\ar\ar\qquad \leq \sum_{s\in J_k(r,t)} \frac{2\Phi_k^{(1)}\big(s, h_k(w_{k}(s))\big)^2}{k u_{k}(s)^2}\le
8{  C_4} k^{-1}\sum_{s\in J_k(r,t)}|\Phi_k^{(1)}\big(s,h_k(w_{k}(s)))|\cr\cr
\ar\ar\qquad \leq 8{  C_4} k^{-1}\big(F_{1,k}(t)-F_{1,k}(r)\big).
\eeqlb
As an application of condition (A) on the sequence $\{h_k(w_{k}(s))\}\subset[0,M_u]$, we have for sufficently large $ k$,
\beqnn
\ar\ar \sum_{s\in J_k(r,t)}u_{k}(s)^{-1}
\Big|\int_{\gamma_k^{-1}(\gamma_k(s-))}
 ^{s-} \phi\big(s,h_k(w_{k}(s))\big)\gamma(\d s)
- \Phi_k^{(1)}\big(s,h_k(w_{k}(s))\big)\Big|\cr\cr
\ar\leq\ar 2\big(F_k(t)-F_k(r)\big).
\eeqnn
Note that $\gamma_k(s-)\ge \gamma_k(r)$ and $\tilde{\lambda}_{k}(s-)
\ge \inf_{s\in [r,t] }v_k(s,t;\lambda)\ge l\ge l_u>0$ for $s\in J_k(r,t)$.
By \eqref{gbk}, \eqref{dphi-M}, {\eqref{w-phi}} and Proposition \ref{prop.Phi-k-M}, for sufficently large $k$,
\beqlb\label{01025}
\ar\ar\sum_{s\in J_k(r,t)}u_{k}(s)^{-1}
\Big| \int_{\gamma_k^{-1}(\gamma_k(s-))}
 ^{s-}\big[ \phi\big(s,h_k(w_{k}(s))\big)
-\phi\big(s,h_k(\tilde{\lambda}_{k}(s-))\big)\big]\gamma(\d s)\Big|
\cr\cr
\ar\ar\qquad \leq 2\beta_k^{-1} {M}'_{\phi}(l_u,M_u) \sum_{s\in J_k(r,t)}
|h_k(w_{k}(s))-h_k(\tilde{\lambda}_{k}(s-))|\cr\cr
\ar\ar\qquad =
2\beta_k^{-1} {M}'_{\phi}(l_u,M_u) \sum_{i\in J_k(r,t)}
|\Phi_k^{(1)}\big(s,h_k(w_{k}(s))\big)|\cr\cr
\ar\ar\qquad \leq 2\beta_k^{-1} {M}'_{\phi}(l_u,M_u)\big(F_{1,k}(t)-F_{1,k}(r)\big),
\eeqlb
since for $s\in J_k(r,t)$,
\beqnn
h_k(w_k(s))=k\big(1-u_k(s)\big)\in[l_u,M_u].
\eeqnn
For sufficently large $k$, $h_k(\tilde{\lambda}_{k}(s-))\ge l/2$ and
\beqnn
\ar\ar\sum_{s\in J_k(r,t)}u_{k}(s)^{-1}
\Big|\int_{\gamma_k^{-1}(\gamma_k(s-))}
 ^{s-}\big[ \phi\big(s,\tilde{\lambda}_{k}(s-)\big)-
 \phi\big(s,h_k(\tilde{\lambda}_{k}(s-))\big)\big]\gamma(\d s)\Big|\cr\cr
 \ar\ar\quad\leq 2{M}'_{\phi}(l/2,U) \sum_{s\in J_k(r,t)}
\int_{\gamma_k^{-1}(\gamma_k(s-))}
 ^{s-}
 |\tilde{\lambda}_{k}(s-)-h_k(\tilde{\lambda}_{k}(s-))|\gamma(\d s)
 \cr\cr
 \ar\ar\quad
 \leq 2{M}'_{\phi}(l/2,U)\frac{U^2}{2k\beta_k} \big(\gamma_k(t)-\gamma_k(r)\big).
\eeqnn
Lemma \ref{l.u} implies for sufficently large $k$,
\beqnn
\ar\ar\sum_{s\in J_k(r,t)}|1-u_{k}(s)^{-1}|
\Big|\int_{\gamma_k^{-1}(\gamma_k(s-))}
 ^{s-} \phi\big(s,\tilde{\lambda}_{k}(s-)\big)\gamma(\d s)
 \Big| \leq \frac{2M_u {M}_{\phi}(U)}{k\beta_k}\big(\gamma_k(t)-\gamma_k(r)\big).
 \eeqnn
Then the desired result can be obtained.
\qed

We conclude this section with an another decomposition of $v_k(r,t;\lambda)$. For $  0\le r\le t$, \eqref{eq.v-de} can
be rewritten as the following equation:
\beqlb\label{vkde2}
v_k(r,t;\lambda) = \lambda - \int_r^t\phi\big(s,v_k(s,t,\lambda)\big)\gamma(\d s) + \sum_{i=1}^4 I^{(i)}_k (r,t;\lambda),
\eeqlb
where
\beqlb
\ar\ar I^{(1)}_k (r,t;\lambda)
= \int_{\gamma_k^{-1}(\gamma_k(t))}^t \phi\big(s,v_k(s,t,\lambda)\big) \gamma(\d s)
- \int_{\gamma_k^{-1}(\gamma_k(r))}^r \phi\big(s,v_k(s,t,\lambda)\big) \gamma(\d s);\nnm\\
 \ar\ar
 I^{(2)}_k (r,t;\lambda) = \sum_{
s\in {J}^+_k(r,t)
}\big[\phi(s,\tilde{\lambda}_{k}(s))\Delta\gamma(s) -\phi^{(2)}_k\big(s, \tilde{\lambda}_{k}(s)\big) \big];\label{I2}\\
\ar\ar
 I^{(3)}_k (r,t;\lambda)= \sum_{s\in J_k(r,t)}\Big[\int_{\gamma_k^{-1}(\gamma_k(s-))}^{s-} \phi(s,v_k(s,t;\lambda))\gamma(\d s)
 - L_{k}(s,t;\lambda) \Big];\label{I3}\\
\ar\ar I^{(4)}_k (r,t;\lambda) = \sum_{
s\in J_k(r,t): \Delta\gamma_k(s)=1} \phi(s,\tilde{\lambda}_{k}(s))\Delta\gamma(s).\nnm
\eeqlb
Denote $I^{(i)}_k (r,t;\lambda)$  by $I^{(i)}_k$ without confusion for $i=1,2,3,4$.  Then under the conditions (A) and (B),
 {  for $\lambda\in[a,b]$ and sufficiently large $k$, we have}
\beqnn
|I^{(1)}_k |
 \ar\le\ar  2\beta_k^{-1}{M}_\phi(U);\\
|I^{(2)}_k |
\ar\leq \ar F_{2,k}(t)-F_{2,k}(r);\\
|I^{(3)}_k |\ar \le \ar
F_{3,k}(t)-F_{3,k}(r),
\eeqnn
and by Lemma \ref{l.cc} and  Lemma \ref{l.cj},
\beqlb\label{I4}
|I^{(4)}_k | \ar \le \ar \sum_{
s\in J_k(r,t): \Delta\gamma_k(s)=1} | \phi(s,\tilde{\lambda}_{k}(s))\Delta\gamma(s)|
 \le F_{4,k}(t)-F_{4,k}(r),\eeqlb
 where
 \beqlb\label{F4}
 F_{4,k}(t):={M}_\phi(U)\sum_{
s\in (0,t]} \Delta\gamma(s)\mbf{1}_{\{s\in (0,t]: \Delta\gamma(s)\le 2\beta_k^{-1}\}}
\eeqlb
 tends to $0$  uniformly for $0\leq t\leq T$ as $k\rightarrow\infty$ by dominated convergence theorem.

\section{Proofs of the main theorems}
\setcounter{equation}{0}
In this section, we will proof Theorem \ref{prop.vk-M} at first. Based on it, we will prove Theorem \ref{t.vk.v} and Theorem \ref{t1.4}.

{\noindent{\it Proof of Theorem \ref{prop.vk-M}.~~}} {  Let $\lambda\in[a,b]$ in the following proof,
and set
\beqnn
\ar\ar U=(b+C_0\gamma(T) +1)\e^{C_0\gamma(T)};\cr\cr
\ar\ar l= 2^{-1}a\prod_{s\in(0,T]} [1+(0\land\Delta\alpha(s))] \exp\{-\|\alpha\|(T)\}>0,
\eeqnn}
where $\alpha$ is defined in \eqref{eq2.20}. Let
$$  t_k:=t_k(t)=\sup\{0\leq s< t:\exists \lambda\in [a,b],~v_k(s,t;\lambda)>U~\mbox{or}~v_k(s,t;\lambda)<l\}$$
and $\sup \emptyset = 0$ by convention. Note that $r\in[0,t]\mapsto v_k(r,t;\lambda)$ is a c\`{a}dl\`{a}g piecewise
function. Then $v_k(t_k,t;\lambda)\in [l,U]$ and $t_k=\gamma_k^{-1}(\gamma_k(t_k))=\gamma_k^{-1}(\gamma_k(t_k-)+1)$.
On the other hand, for $s\in[\gamma^{-1}_k(\gamma_k(t_k-)),t_k)$, $v_k(s,t;\lambda)= v_k(\gamma^{-1}_k(\gamma_k(t_k-)),t;\lambda) = v_k(t_k-,t;\lambda)$.
Then for $r_0\in [\gamma^{-1}_k(\gamma_k(t_k-)),t_k)$,
\beqlb\label{vtkde}
v_k(t_k-,t;\lambda)\ar =\ar v_k(r_0,t;\lambda) \cr
\ar=\ar v_k(t_k,t;\lambda) - \sum_{s\in J_k(r_0,t_k)}L_{k}(s,t;\lambda) - \sum_{s\in {J}^+_k(r_0,t_k)} \phi^{(2)}_k(s,v_k(s,t;\lambda)).
\eeqlb
By  the proof of Lemma \ref{l.u} and the boundedness of $v_k(t_k,t;\lambda)$,
\beqnn
u_{k}(t_k)=\left\{
\begin{array}{ll}
\exp\{-k^{-1}[v_k(t_k,t;\lambda)-\phi^{(2)}(t_k,v_k(t_k,t;\lambda))]\}, & \mbox{if~$\Delta\gamma_k(t_k)>1$;} \\
\exp\{-k^{-1}v_k(t_k,t;\lambda)\} , & \mbox{if~$\Delta\gamma_k(t_k)=1$.}
\end{array}
\right.
\eeqnn
If $\Delta\gamma_k(t_k)>1$, then $t_k\in {J}^+_k(r_0,t)$ and by a similar argument in Lemma \ref{l.u}, for sufficiently large $  k\ge N_3$,
\beqnn
{  3 l_u/2}\le v_k(t_k,t;\lambda)-\phi^{(2)}(t_k,v_k(t_k,t;\lambda)) \le M_u,
\eeqnn
where  {  $M_u$ and $l_u$ are defined} in \eqref{Mulul'}. Hence for sufficiently large $k$,
$|k(1-u_{k}(t_k)|\le M_u$ and $u_{k}(t_k)\ge 1/2$. Then by the proof Proposition \ref{prop.Phi-k-M}, {  there
exists $ F_{1,k}$} such that, for sufficiently large $k$,
\beqnn
\Big|\frac{\Phi_k^{(1)}\big(  t_k , k(1-u_{k}(t_k))\big)}{ku_{k}(t_k)}\Big|
\le  \frac{2F_{1,k}(T)+2{M}_{\phi}(M_u)/\beta_k}{k}\le { \frac{l_u}{8k}\le \frac{1}{2}},
\eeqnn
and
\beqnn
\ar\ar \Big| \sum_{s\in J_k(r_0,t_k)}L_{k}(s,t;\lambda)\Big|
=  | L_{k}(t_k,t;\lambda)| = \Big|k\log\Big(\frac{\Phi_k^{(1)}\big( t_k, k(1-u_{k}(t_k))\big)}{ku_{k}(t_k)}+1\Big)\Big|\cr
\ar\ar\qquad \le
\Big| \frac{\Phi_k^{(1)}\big( t_k, k(1-u_{k}(t_k))\big)}{u_{k}(t_k)}\Big|
+ 2k\Big|\frac{\Phi_k^{(1)}\big(  t_k , k(1-u_{k}(t_k))\big)}{ku_{k}(t_k)}\Big|^2 \\
\ar\ar\qquad \le  \frac{l_u}{8 }+  \frac{l_u}{8 } = \frac{l_u}{4}.
\eeqnn
Let $ U'=U+M_\phi(U)\gamma(T)+l$. Then by an elementary calculation, there exists $  K_0=K_0(T)$ such that for $k\ge K_0$, we have
\beqlb\label{tk-}
  v_k(t_k-,t;\lambda)\in [l_u,U'].
\eeqlb If $t_k>0$, then for $r\in[\gamma^{-1}_k(\gamma_k(t_k-)),t]$, by \eqref{vkde2},
\beqlb\label{vkUb}
  v_k(r,t,\lambda) \ar\leq\ar \lambda- \int_{r}^t\phi\big(s,v_k(s,t,\lambda)\big)\gamma(\d s)+\sum_{i=1}^4 |I_k^{(i)}(r,t;\lambda)|\cr
\ar\le\ar
\lambda + C_0\gamma(T) + C_0\int_{r}^t v_k(s,t;\lambda)\gamma(\d s)  + \sum_{i=1}^4|I^{(i)}_k(r,t;\lambda)|,
\eeqlb
since $\phi(s,\lambda)\geq -|b_1(s)|\lambda-\int_1^{\infty}m(s,\d z)\geq -(C_0+C_0\lambda)$. By the discussions in the previous section
and the boundedness of $v_k(t_k-,t;\lambda)$ and $v_k(t_k,t;\lambda)$, there exists $  K_1=K_1(T,l_u,U')$ such
that $\sum_{i=1}^4 |I^{(k)}_i(r,t;\lambda)|\leq 1$ for $r\in[\gamma^{-1}_k(\gamma_k(t_k-)),t]$ and $k\ge K_1$. By applying Gronwall's
inequality to \eqref{vkUb} for $k\ge K_1$ and $r\in[\gamma^{-1}_k(\gamma_k(t_k-)),t]$, we get
\beqlb\label{u}
 v_k(r,t;\lambda)\leq [\lambda+C_0\gamma(T)+1]\exp\{C_0(\gamma(t)-\gamma(r))\}\leq U.
 \eeqlb
On the other hand, let $r\mapsto\pi_k(r,t;\lambda)$ be the solution of
\beqnn
\pi_k(r,t;\lambda) = \lambda + \sum_{i=1}^4I^{(i)}_k(r,t;\lambda) + \int_r^t\pi_k(s,t;\lambda)\alpha(\d s),
\qquad r \in[\gamma^{-1}_k(\gamma_k(t_k-)),t].
\eeqnn
Following a similar proof of Proposition 2.3 in Fang and Li (2022+), we have $ v_k(r,t;\lambda)\ge \pi_k(r,t;\lambda)$.
By Proposition 2.1 in Fang and Li (2022+),
\beqlb\label{0603-1}
\pi_k(r,t;\lambda) = \e^{\zeta(t)-\zeta(r)} \lambda - \int_r^t\e^{\zeta(s-)-\zeta(r)}G(\d s),
\eeqlb
where $\zeta_c(t)= \alpha_c(t)$, $\Delta\zeta(t)= \log [1+\Delta\alpha(t)]$ and $G(r)=\lambda + \sum_{i=1}^4I^{(i)}_k(r,t;\lambda)$. Then
\beqlb\label{0603}
\Big|\int_r^t\e^{\zeta(s-)-\zeta(r)}G(\d s)\Big| \le  \sum_{i=1}^4\Big|\int_r^t\e^{\zeta(s-)-\zeta(r)}I^{(i)}_k(\d s,t;\lambda)\Big|.
\eeqlb
By \eqref{Lambda}, $\Lambda$ is locally bounded by $\Lambda_T$, hence is bounded for sufficiently large $k$ such that $\beta_k^{-1}\Lambda_T<1/2$.
Moreover,  for $s\in J_k(r,t)$ and $\gamma_k^{-1}(\gamma_k(s-))\le u< s$, since $|\Delta\alpha(u)|\leq\Lambda_T\Delta \gamma(u)\leq \Lambda_T/\beta_k\leq 1/2$
for sufficiently large $k$, we have for $\Delta\alpha(u)<0$,
\beqlb\label{log}
|\log\big(1+\Delta\alpha(u)\big)|\ar=\ar\log \Big(1-\frac{\Delta\alpha(u)}{1+\Delta\alpha(u)}\Big)\cr\cr
\ar\leq\ar \log\big(1-2\Delta\alpha(u)\big)\leq|2\Delta\alpha(u)|.
\eeqlb
For $\Delta\alpha(u)\geq 0$, it's obvious that
\beqnn
\log\big(1+\Delta\alpha(u)\big)\leq \Delta\alpha(u).
\eeqnn
Combining it with \eqref{log}, we have
\beqnn
|\log\big(1+\Delta\alpha(u)\big)|\leq|2\Delta\alpha(u)|.
\eeqnn
Consequently, for $s\in J_k(r,t)$ and $\gamma_k^{-1}(\gamma_k(s-))\le u< s$,
 \beqlb\label{zeta}
 |\zeta(s-)-\zeta(u-)|\ar\leq\ar
 |\alpha_c(s-)-\alpha_c(u-)|+\sum_{u\in [\gamma_k^{-1}(\gamma_k(s-)),s)}|\log\big(1+\Delta\alpha(u)\big)|\cr\cr
 \ar\leq\ar \Lambda_T /\beta_k+2\sum_{u\in [\gamma_k^{-1}(\gamma_k(s-)),s)}|\Delta\alpha(u)|\cr\cr
 \ar\leq\ar
  \Lambda_T /\beta_k+2\Lambda_T\sum_{u\in [\gamma_k^{-1}(\gamma_k(s-)),s)}|\Delta\gamma(u)|\cr\cr
  \ar\leq\ar
  3 \Lambda_T/\beta_k.
 \eeqlb
 Since
 \beqlb\label{deltagamma}
\gamma_k(t)/\beta_k\leq \gamma(\gamma_k^{-1}(\gamma_k(t)))\leq \gamma(t)<\big(\gamma_k(t)+1\big)/\beta_k,
\eeqlb
we have
 \beqlb\label{gamma}
 \int_{\gamma_k^{-1}(\gamma_k(r))}^r \gamma(\d s)\leq 1/\beta_k;\quad
  \int_{\gamma_k^{-1}(\gamma_k(t))}^t \gamma(\d s)\leq 1/\beta_k.
  \eeqlb
  For $r \in[\gamma^{-1}_k(\gamma_k(t_k-)),t]$ and $r<s\leq t$,
 \beqnn
 \e^{\zeta(s-)-\zeta(r)}
 \ar\leq\ar
 \e^{\zeta_c(s-)-\zeta_c(r)}\prod_{r\leq u<s}\e^{\Delta\zeta(u)}\le\e^{\|\alpha\|(T)}\prod_{s\in(0,T]}(1+|\Delta\alpha(s)|).
 \eeqnn
 On the other hand, by \eqref{gamma},
 \beqnn
 \sum_{\gamma_k^{-1}(\gamma_k(r))<u\leq r}\Delta\gamma(u)\leq 1/\beta_k.
 \eeqnn
 Then by a similar argument with \eqref{zeta},  for $\gamma_k^{-1}(\gamma_k(r))<s\leq r$, we have
 \beqnn
  \e^{\zeta(s-)-\zeta(r)}\leq\e^{3 \Lambda_T/\beta_k},
  \eeqnn
  for sufficiently large $k$.
  Consequently, for sufficiently large $k$ and $\gamma_k^{-1}(\gamma_k(r))<s\leq t$,
  \beqlb\label{zeta2}
  \e^{\zeta(s-)-\zeta(r)}\leq \e^{3 \Lambda_T/\beta_k}\vee\big[\e^{\|\alpha\|(T)}\prod_{s\in(0,T]}(1+|\Delta\alpha(s)|)\big]
  \leq 2\e^{\|\alpha\|(T)}\prod_{s\in(0,T]}(1+|\Delta\alpha(s)|).
  \eeqlb
  Note that by \eqref{u}, for $\gamma_k^{-1}(\gamma_k(r))<s\leq t$,
 \beqnn
 |\phi\big(s,v_k(s,t;\lambda)\big)|\leq M_\phi(U).
 \eeqnn
 Since
 \beqnn
 I_k^{(1)}(\d s,t;\lambda)=-\phi\big(s,v_k(s,t;\lambda)\big)\gamma(\d s)
 +\mbf{1}_{\{\Delta\gamma_k(s)\geq 1\}}\int_{\gamma_k^{-1}(\gamma_k(s-))}^{\gamma_k^{-1}(\gamma_k(s))}\phi\big(u,v_k(u,t;\lambda)\big)\gamma(\d u),
\eeqnn
 it follows from \eqref{zeta2} that for sufficiently large $k$,
\beqnn
\ar\ar\Big|\int_r^t\e^{\zeta(s-)-\zeta(r)}I^{(1)}_k(\d s,t;\lambda)\Big|\cr
\ar\ar=
\big|-\int_r^t\e^{\zeta(s-)-\zeta(r)}\phi\big(s,v_k(s,t;\lambda)\big)\gamma(\d s)\cr
\ar\ar\quad
+\sum_{s\in J_k(r,t)}\int_{\gamma_k^{-1}(\gamma_k(s-))}^{\gamma_k^{-1}(\gamma_k(s))}\e^{\zeta(s-)-\zeta(r)}
\phi\big(u,v_k(u,t;\lambda)\big)\gamma(\d u)\big|\cr
\ar\ar\leq
\big|\int_{\gamma_k^{-1}(\gamma_k(r))}^{\gamma_k^{-1}(\gamma_k(t))}\e^{\zeta(s-)-\zeta(r)}\phi\big(s,v_k(s,t;\lambda)\big)
\gamma(\d s)-\int_r^t\e^{\zeta(s-)-\zeta(r)}\phi\big(s,v_k(s,t;\lambda)\big)\gamma(\d s)\big|\cr
\ar\ar\quad
+ \sum_{s\in J_k(r,t)}\Big|  \int_{\gamma_k^{-1}(\gamma_k(s-))}^{s} (\e^{\zeta(u-)-\zeta(r)}-\e^{\zeta(s-)-\zeta(r)})\phi\big(u,v_k(u,t;\lambda)\big)
 \gamma(\d u)\Big|\cr
\ar\le \ar
\Big|\int_{\gamma_k^{-1}(\gamma_k(r))}^r \e^{\zeta(s-)-\zeta(r)}\phi\big(s,v_k(s,t;\lambda)\big) \gamma(\d s)\Big| \cr
\ar\ar\qquad
+ \Big|\int_{\gamma_k^{-1}(\gamma_k(t))}^t \e^{\zeta(s-)-\zeta(r)} \phi\big(s,v_k(s,t;\lambda)\big) \gamma(\d s)\Big| \cr
\ar\ar\qquad  + \sum_{s\in J_k(r,t)}\Big|  \int_{\gamma_k^{-1}(\gamma_k(s-))}^{s} \e^{\zeta(u-)-\zeta(r)}(1-\e^{\zeta(s-)-\zeta(u-)})\phi\big(u,v_k(u,t;\lambda)\big)
 \gamma(\d u)\Big|\cr
\ar\le \ar 4\beta_k^{-1}M_\phi(U)\e^{\|\alpha\|(T)}\prod_{s\in(0,T]}(1+|\Delta\alpha(s)|)\cr
\ar\ar\qquad
+ 2M_\phi(U)\e^{\|\alpha\|(T)}\prod_{s\in(0,T]}(1+|\Delta\alpha(s)|)
\sum_{s\in J_k(r,t)}   \int_{\gamma_k^{-1}(\gamma_k(s-))}^{s} |1-\e^{\zeta(s-)-\zeta(u-)}| \gamma(\d u)\cr
\ar\le \ar
4\beta_k^{-1}M_\phi(U)\e^{\|\alpha\|(T)}\prod_{s\in(0,T]}(1+|\Delta\alpha(s)|)\cr
\ar\ar\qquad + 2M_\phi(U)\e^{\|\alpha\|(T)}\prod_{s\in(0,T]}(1+|\Delta\alpha(s)|) (\e^{3\beta_k^{-1}\Lambda_T}-1)\gamma(T),
\eeqnn
where the last inequality follows from $|\e^x-1|\leq\e^{|x|}-1$ for $x\in\mbb{R}$.
On the other hand, combining \eqref{I2}, \eqref{I3} and \eqref{I4}, by \eqref{tk-}, Lemma \ref{l.cj} and Lemma \ref{l.cc},
 for sufficiently large $k$, we have
\beqnn
\ar\ar\Big|\sum_{i=2}^4\int_r^t\e^{\zeta(s-)-\zeta(r)}I^{(i)}_k(\d s,t;\lambda)\Big|\cr\cr
\ar\ar\leq 2\e^{\|\alpha\|(T)}\prod_{s\in(0,T]}(1+|\Delta\alpha(s)|)
\sum_{i=2}^4\int_r^t\|I^{(i)}_k\|(\d s,t;\lambda)\cr\cr
\ar\ar\leq 2\e^{\|\alpha\|(T)}\prod_{s\in(0,T]}(1+|\Delta\alpha(s)|)\sum_{i=2}^4[F_{i,k}(t)-F_{i,k}(r)],
\eeqnn
which tends to $0$ as $k\rightarrow\infty$. Therefore $\big|\int_r^t\e^{\zeta(s-)-\zeta(r)}G(\d s)\big|$ can be
sufficiently small for sufficiently large $k$ and $r \in[\gamma^{-1}_k(\gamma_k(t_k-)),t]$  by \eqref{0603}. Hence there
exists $  N_1=N_1(T,l_u,U')\ge K_1$ such that for $r \in[\gamma^{-1}_k(\gamma_k(t_k-)),t]$
{  and $k\ge N_1$}, $\big|\int_r^t\e^{\zeta(s-)-\zeta(r)}G(\d s)\big|\leq l$ and then combining with \eqref{u} and
observing that $\e^{\zeta(t)-\zeta(r)} \lambda\ge 2l$, we have
 $$ U\ge v_k(r,t;\lambda)\ge \pi_k(r,t;\lambda)\ge 2l-l=l,$$
 which contradicts   the definition of $t_k$ and implies that $t_k=0$. Therefore we obtain the result {  with  $ k\geq N_1$}. \qed

{\noindent{\it Proof of Theorem \ref{t.vk.v}.~~}}
   By discussion in the last section above,   $I^{(1)}_k,~I^{(2)}_k,~I^{(3)}_k,~I^{(4)}_k$ tend
to $0$ uniformly {  for $0\le r\le t\le T$ and $\lambda\in[a,b]$} as $k\rightarrow\infty$.
Consequently, for sufficiently large $k$ depending on $  a$ and $  b$, a standard argument leads to
 { \beqnn
 \ar\ar\sup_{m,n\geq k\atop \lambda\in[a,b]}|v_m(r,t;\lambda)-v_n(r,t;\lambda)|\cr\cr
 \ar\ar\quad\leq {M}'_{\phi}(l,U)\int_r^t\sup_{m,n\geq k\atop\lambda\in[a,b]}|v_m(s,t;\lambda)-v_n(s,t;\lambda)|\gamma(\d s)+ \epsilon_{m,n}(a,b),
 \eeqnn}
where $ \epsilon_{m,n}(a,b)$ converges to $0$ as $k\rightarrow\infty$.
As an application of Gronwall's inequality,  we find that
$v_n(s,t;\lambda)$ converges to some limit uniformly on $  (s,\lambda)\in [0,t]\times[a,b]$ as $n\rightarrow\infty$. Taking limits of \eqref{vkde2},
we can complete the proof with the uniqueness of the solution \eqref{v.fll}.\qed

{\noindent{\it Proof of Theorem \ref{t1.4}.~~}} By Theorem \ref{t.vk.v}, it's sufficient to prove the tightness of $\{X_k\}$.
 {  Let $ T\geq 0$ and $0\leq r\leq t\leq T$. Then by} \eqref{eq.v-de}, we have
\beqlb\label{vlambda}
|\lambda-v_k(r,t;\lambda)|\ar=\ar|\sum_{s\in J_k(r,t)} L_{k}(s,t,\lambda) - \sum_{s\in {J}^+_k(r,t)}
 \phi^{(2)}_k\big(s,v_k(s,t;\lambda)\big)|\cr\cr
 \ar\leq\ar |I^{(2)}_k| +|I^{(3)}_k| +|I^{(4)}_k|+
 \Big|\int_{\gamma_k^{-1}(\gamma_k(r))}^{\gamma_k^{-1}(\gamma_k(t))} \phi\big(s,v_k(s,t,\lambda)\big) \gamma(\d s)\Big|\cr\cr
 \ar\leq\ar F_{2,k}(t)-F_{2,k}(r)+F_{3,k}(t)-F_{3,k}(r)\cr\cr
 \ar\ar\qquad+F_{4,k}(t)-F_{4,k}(r)+F_{5,k}(t)-F_{5,k}(r),
 \eeqlb
where  $F_{4,k}$ is given by \eqref{F4} and $F_{5,k}$ is an increasing c\`adl\`ag function defined by
\beqnn
 F_{5,k}(t)={M}_{\phi}(U)\gamma(\gamma_k^{-1}(\gamma_k(t))).
\eeqnn
Moreover $F_{5,k}(t)$ converges to ${M}_{\phi}(U)\gamma(t)$ uniformly for $t\in[0,T]$ as $k\rightarrow\infty$ by \eqref{deltagamma}.
 By Lemma \ref{l.cj} and Lemma \ref{l.cc}, $F_{2,k}(t)$ and $F_{3,k}(t)$ converge to $0$ uniformly {  for $  t\in[0, T]$}
 as $k\rightarrow\infty$ indeed. Consequently,  there exists a sequence of increasing c\`adl\`ag
 functions $F_{6,k}$ such that for sufficiently large $k$,
\beqlb\label{vtight}
 |\lambda-v_k(r,t;\lambda)|\leq F_{6,k}(t)-F_{6,k}(r),\ \quad 0\leq r\leq t\leq T,
 \eeqlb
 where $F_{6,k}(t)$ converges to ${M}_{\phi}(U)\gamma(t)$ uniformly {  for $  t\in [0,T]$} as $k\rightarrow\infty$. Hence
 $F_{6,k}$ converges to ${M}_{\phi}(U)\gamma$ {  in $ D[0,\infty)$} with the Skorohod topology as $k\rightarrow\infty$.
  By the proof of claim IV in Bansaye and Simatos (2015), for sufficiently large $k$,
   \beqnn
   \mbf{E}[d^2\Big(X_k(t),X_k(r)\Big)|X_k(r)=x_0]
   \ar=\ar
   \e^{-2x_0}[2(1-\e^{x_0(1-v_k(r,t;1))})-(1-\e^{x_0(2-v_k(r,t;2))})]\cr
   \ar\leq\ar (2x_0\e^{-2x_0}+1)|1-v_k(r,t;1)|\cr\ar\ar\qquad +(x_0\e^{-2x_0}+1)|2-v_k(r,t;2)|\cr
   \ar\leq\ar 2(1+\e^{-1})(F_{6,k}(t)-F_{6,k}(r)).
   \eeqnn
  By Theorem 1 of Bansaye, Kurtz and Simatos (2016), we obtain {  the tightness } since $[0,\infty]$ endowed with the metric $d$ is compact.
  {  On the other hand, since for $0\leq s\leq t$ and $\lambda,\mu>0$,
  \beqnn
  \lim_{k\rightarrow\infty}\mbf{E}[\e^{-\lambda X_k(t)-\mu X_k(s)}]
  \ar=\ar
  \lim_{k\rightarrow\infty}\mbf{E}\big[\mbf{E}[\e^{-\lambda X_k(t)-\mu X_k(s)}|X_k(0)]\big]\cr\cr
  \ar=\ar
  \lim_{k\rightarrow\infty}\mbf{E}\big[\mathrm{exp}\bigl\{-X_k(0) v_k\big(0,s;\mu+v_k(s,t;\lambda)\big)\bigr\}\big]\cr\cr
  \ar=\ar \mbf{E}\big[\mathrm{exp}\bigl\{-X(0) v\big(0,s;\mu+v(s,t;\lambda)\big)\bigr\}\big]=\mbf{E}[\e^{-\lambda X(t)-\mu X(s)}],
  \eeqnn
by  induction, one can show that the finite-dimensional distributions of $\{X_k(s):s\geq 0\}$ converge to that of $\{X(s):s\geq 0\}$. Hence we obtain the weak convergence.\qed

\bigskip

\textbf{Acknowledgements.} The authors want to thank Professor V.A. Vatutin
for his careful reading of the paper and a number of helpful corrections.}

\end{document}